\documentclass{article}
\usepackage{amssymb}
\usepackage{amsmath}
\newtheorem{theorem}{Theorem}[section]
\newtheorem{cor}[theorem]{Corollary}

\newtheorem{lem}[theorem]{Lemma}
\newtheorem{pro}[theorem]{Proposition}

\newtheorem{rem}[theorem]{Remark}

\def\L2{L^2([ 0,2\pi))}
\def\L2R{L^2(\bold{R})}

\def\cA{{\cal A }}
\def\cB{{\cal B }}

\def\cP{{\cal P }}

\def\CC{{\Bbb C}}
\def\NN{{\Bbb N}}
\def\QQ{{\Bbb Q}}
\def\RR{{\Bbb R}}
\def\TT{{\Bbb T}}

\def\ZZ{{\Bbb Z}}

\begin{document} 
\title{Fixed point and spectral characterization of finite dimensional
$C^*$--algebras\thanks{Supported by the Thailand Research Fund
under grant BRG50800016. The second author was supported by the
Development and Promotion of Science and Technology Talent Project
(DPST).}}

\author{S. Dhompongsa$^1$\thanks{Corresponding author.\newline  E-mail addresses :
sompongd@chiangmai.ac.th (S. Dhompongsa), g4865050@cm.edu (W.
Fupinwong), matwml@nus.edu.sg (W. Lawton). \newline Key words and
phrases: fixed point property for nonexpansive maps, spectrum,
minimal projection, Gelfand-Naimark-Segal theorem.
}{,~~W. Fupinwong$^1$ and W. Lawton$^2$ }\\{\small
$^1$Department of Mathematics, Faculty of Science, Chiang Mai
University,}\\{\small Chiang Mai 50200, Thailand}\\{\small
$^2$Department of Mathematics, National University of Singapore,
}\\{\small 2 Science Drive 2, Singapore 117543}}

\date{}

\maketitle
\begin{abstract}
We show that the following conditions on a $C^*$--algebra are
equivalent: (i) it has the fixed point property for nonexpansive
mappings, (ii) the spectrum of every self adjoint element is
finite, (iii) it is finite dimensional. We prove that (i) implies
(ii) using constructions given by Goebel, that (ii) implies (iii)
using projection operator properties derived from the spectral and
Gelfand-Naimark-Segal theorems, and observe that (iii) implies (i)
by Brouwer's fixed point theorem.
\end{abstract}

\maketitle

\section{Statement of Result}
\label{introduction}
\setcounter{equation}{0}
We let $\CC,$ $\RR,$ $\QQ,$ $\ZZ,$ $\NN = \{1,2,3,4,...\},$ and
$\TT = \{ \, z \in \CC \, : \, |z| = 1 \, \}$ denote the fields of
complex, real, and rational numbers, the ring of integers, the set
of natural numbers, and the circle group. If $\cB$ is a Banach
space and $K \subset \cB$ then a map $M \, : \, K \rightarrow K$
is {\it nonexpansive} if $||M(x) - M(y)|| \leq ||x-y||, \, x, y
\in K,$ and $\cB$ has the {\it fixed point property for
nonexpansive maps} (FPP) if $M$ has a fixed point whenever $M$ is
nonexpansive and $K$ is a nonempty bounded closed convex subset of
$\cB,$ \cite{goebel}. Throughout this paper $\cA$ denotes a (not
necessarily unital) $C^*$--algebra with $\cA \neq \{ 0\}$ and
$\widetilde \cA$ denotes $\cA$ if $\cA$ is unital and the
unitization of $\cA$ otherwise and $I$ denotes its unit.
$\text{Inv}(\widetilde \cA)$ denotes the set of invertible
elements in $\widetilde \cA$ and the {\it spectrum} of an element
$A \in \cA$ is defined by
    $\hbox{spec}(A) = \{ \, \mu \in \CC \, : \, \mu \, I - A \notin
    \text{Inv}(\widetilde \cA)\, \},$
\cite{murphy}. We let $\cA_s, \cA_u,$ and $\cA_n$ denote the sets
of self-adjoint, unitary, and normal elements in $\cA.$ The
following result relates fixed point and spectral properties of
$C^*$--algebras.

\bigskip
\begin{theorem}
\label{pro:GNS} The following properties of a $C^*$--algebra $\cA$
are equivalent:
\begin{enumerate}
\item $\cA$ has the FPP,
\item ${spec}\, (A)$ is finite whenever $A \in \cA_s,$
\item $\cA$ has finite dimension.
\end{enumerate}
\end{theorem}
\textbf{Proof.} Proposition \ref{pro1} shows that property (1)
implies property (2), Proposition \ref{pro2} shows that property
(2) implies property (3), and Brouwer's Fixed Point Theorem
(\cite{goebel}, Theorem 7.1) shows that property (3) implies
property (1).
\hfill$\square$\\
\section{Derivation of Proposition \ref{pro1}}
\label{derivations1} \setcounter{equation}{0}
Throughout this section $S$ denotes a compact subset of $\CC.$
We observe that set $C(S),$ consisting of continuous complex valued functions on $S$ and equipped with adjoint $A^* = \overline A$ and norm $||A|| = \max_{\, s \in S} |A(s)|,$ is a unital abelian $C^*$--algebra and that for every $p \in S$ the subset $I_p(S),$ consisting of functions in $C(S)$ that vanish at $p,$ is an ideal of $C(S).$
If $A \in \cA$ then $\hbox{spec}(A)$ is a compact subset of the
disc of radius $||A||$ centered at $0.$ Furthermore, $\cA_s$ and
$\cA_u$ are subsets of $\cA_n,$ $A \in \cA_s$ implies
$\hbox{spec}(A) \subseteq [ \, -||A||,||A|| \, ],$ and $A \in
\cA_u$ implies $\hbox{spec}\, (A) \subseteq \TT.$ For $A \in \cA,$
the $C^*$--subalgebra $C^*(A)$ of $\cA$ generated by $A$ and the
$C^*$--subalgebra $C^*(A,I)$ of $\widetilde \cA$ generated by $A$
and $I$ are abelian iff $A \in \cA_n.$
\begin{lem}
\label{abelian}
If $A \in \cA_n$ and
$z \, : \, \hbox{spec} \, (A) \rightarrow \CC$
denotes the inclusion map then there exists an isometric $*$--isomorphism
$\phi \, : \, C^*(A,I) \rightarrow C(\hbox{spec} \, (A))$
such that $\phi(A) = z$ and $\phi(I) = 1.$ If $A \in
\cA_s$ then
$0 \notin \hbox{spec} \, (A) \implies \phi(C^*(A)) = C(\hbox{spec} \, (A))$
and
$0 \in \hbox{spec} \, (A) \implies \phi(C^*(A)) = I_0(\hbox{spec} \, (A)).$
\end{lem}
\textbf{Proof.}
The Gelfand theorem for abelian $C^*$--algebras implies that the maximal ideal space $\Omega(C^*(A,I))$ of $C^*(A,I),$ identified with the set of multiplicative linear functionals $\omega \, : \, C^*(A,I) \rightarrow \CC$ equipped with the weak$^*$ topology, is a compact Hausdorff space and that $C^*(A,I)$ is isomorphic to $C(\Omega(C^*(A,I)))$ under the Gelfand transform $B \rightarrow \widehat B$ defined by $\widehat {B}(\omega) = \omega(B).$ We observe that $\widehat {A} \, : \, \Omega(C^*(A,I)) \rightarrow \hbox{spec} \, (A)$ is a homeomorphism and that the function $\phi \, : \, C^*(A,I) \rightarrow C(\hbox{spec} \, (A))$ defined by  $\phi(B) = \widehat {B} \circ (\widehat {A})^{-1}, \, B \in C^*(A,I),$ satisfies the asserted properties, see also (\cite{murphy}, Theorem 2.1.13), (\cite{strocchi}, Proposition 1.5.3). The second assertion follows from the Weierstrass approximation theorem.
\hfill$\square$\\

\bigskip
The following result generalizes Goebel's Example 2.3 in
\cite{goebel}.
\begin{lem}
\label{goebel1}
If $S$ is infinite then $C(S)$ does not
have the FPP.
\end{lem}
\textbf{Proof.} Choose a limit point $p \in S$ and define the set
$$
    K = \{ \, A \in C(S) \, : \,
    \ 0 \leq \|A\| \leq 1, \ A(p) = 1 \, \}.
$$
Choose $B \in K$ with $|B(s)| < 1$ for $s \neq p$ and define the
map
$M_B \, : \, K \rightarrow K$
by
$ M_B(x) =  BA.$
Then $K$ is a bounded closed convex set and the map $M_B$ is
nonexpansive and does not have a fixed point. \hfill$\square$\
\begin{lem}
\label{goebel2}
If $S$ is infinite and $p \in S$ then $I_p(S)$ does not have the FPP.
\end{lem}
\textbf{Proof.} If there exists a limit point $q \in S$ such that
$q \neq p$ then an argument similar to that used in Lemma
\ref{goebel1} shows that $I_p(S)$ does not have the FPP. Otherwise
$p$ is the only limit point of $S$ and then $I_p(S)$ is
isomorphic, as a Banach space, to the Banach space $c_0$
consisting of all complex valued convergent sequences. Let $K$ be
the unit ball in $c_0.$ Then $K$ is a bounded close convex set and
Goebel's Example 4.1 in \cite{goebel} shows that the mapping
$M \, : \, K \rightarrow K$
defined by
$M(A) = (1-||a||,a_1,a_2,...), \ a \in c_0$
is nonexpansive and does not have a fixed point. Therefore $c_0,$
and hence $I_p(S),$ does not have the FPP. \hfill$\square$\
\begin{pro}
\label{pro1}
If $\cA$ has the FPP and $A \in \cA$ then $\hbox{spec}\, (A)$ is
finite.
\end{pro}
\textbf{Proof.} Assume that $A \in \cA_s$ and $\hbox{spec}\, (A)$
is infinite. It suffices to show that $\cA$ does not have the FPP.
Lemma \ref{abelian} shows that $C^*(A)$ is either isomorphic to
$C(\hbox{spec}\, (A))$ or isomorphic to $I_0(\hbox{spec}\, (A)).$
Therefore Lemma \ref{goebel1} and Lemma \ref{goebel2} show that
the subalgebra $C^*(A)$ of $\cA,$ and hence $\cA,$ does not have
the FPP. \hfill$\square$\
\section{Projections}
\label{projections} \setcounter{equation}{0}
The set of projections
    $\cP = \{ \, P \in \cA_s \, : \, P^2 = P \, \}$
in $\cA$ admits the partial order
    $P \leq Q  \mbox{ iff } PQ = QP = P.$
We write $P < Q$ if $P \leq Q$ and $P \neq Q.$ We say $P$ and $Q$ are orthogonal if $PQ = QP = 0$ and then write $P \perp Q.$ If $P, Q \in \cP$ then $P + Q \in \cP$ iff $P \perp Q.$ A projection $P$ is maximal if there does not exist a projection $Q$  such that $P < Q.$ If $\cA$ is unital then $I$ is the unique maximal projection. A projection $P \in \cP$ is minimal if $P \neq 0$ and there does not exist a nonzero $Q \in \cP$ such that $Q < P.$ Otherwise $P-Q$ is a nonzero projection and $P = Q + (P-Q).$
\begin{lem}
\label{spec}
If $A \in \cA_s,$ $A \neq 0,$ and $\hbox{spec}\, (A)$ is finite,
then there exist pairwise orthogonal nonzero projections
$P_1,...,P_d \in \cP \cap C^*(A)$
and nonzero elements
$\lambda_1,...,\lambda_d \in \hbox{spec}\, (A)$
such that
    $A = \sum_{j=1}^{d} \lambda_j \, P_j.$
Furthermore, $\overline {AH} = AH$ and there exists a $P_A \in \cP$ such that $P_AH = AH.$
\end{lem}
\textbf{Proof.} By assumption $\hbox{spec}\, (A)$ is finite and
since $A \neq 0$ it contains at least one nonzero element. Let $d
\in \NN,$ let $\lambda_1,...,\lambda_d$ be the distinct nonzero
elements in $\hbox{spec}\, (A),$ and define $\pi_j \in
C(\hbox{spec}\, (A)), j = 1,...,d$ by $\pi_j(\lambda_i) = 1$ if
$i=j$ and zero otherwise. Let $\phi \, : \, C^*(A,I) \rightarrow
C(\hbox{spec} \, (A))$ be the isometric $*$--isomorphism in Lemma
\ref{abelian} and define $P_j = \phi^{-1}(\pi_j), j = 1,...,d.$
Then Lemma \ref{abelian} ensures that $P_1,..,P_d$ are pairwise
orthogonal minimal projection operators in $C^*(A)$ and that their
sum equals $A.$ The second assertion follows by defining $P_A =
\sum_{j=1}^{d} P_j.$ \hfill$\square$\
\begin{rem}
\label{spectral}
Lemma \ref{spec} is a special case of the spectral theorem for normal elements in a $C^*$--algebra \cite{murphy} and for normal operators on a Hilbert space \cite{halmos}. Furthermore, for each $j = 1,...,d,$
    $P_j = AL_j(A)/\lambda_j$ where $L_j$
is the Lagrange interpolating polynomial of degree $d-1$ that satisfies
    $L_j(\lambda_i) = 1$ if $i=j$ and equals zero otherwise.
\end{rem}
\begin{lem}
\label{AN}
If $\hbox{spec}\,  \, (A)$ is finite for every $A \in \cA_s$ then
every sequence in $\cP$ that either strictly increases or
strictly decreases must terminate.
\end{lem}
\textbf{Proof.} Assume to the contrary that there exists an
infinite sequence $P_n, \ n \in \NN$ that is either strictly
increasing or strictly decreasing and construct the sequence of
nonzero projections $Q_n \in \cP$ by $Q_n = P_{n+1} - P_n$ if
$P_n$ increases and $Q_n = P_{n} - P_{n+1}$ if $P_n$ decreases.
Then $A_k = \sum_{n=1}^{k} (1/n)Q_n, \, k \in \NN$ is a Cauchy
sequence in $\cA_s$ that converge to $A \in \cA_s$ and
$\hbox{spec}\, (A) = \{ \, 1/n \, : \, n \in \NN \, \} \cup \{ \,
0 \, \}$ is infinite. This contradiction concludes the proof.
\hfill$\square$\
\begin{cor}
\label{maxmin}
If $\hbox{spec}\, (A)$ is finite for every $A \in \cA_s$ then
for every nonzero $P \in \cP$ there exists a maximal $P_{max} \in \cP$ and
a minimal $P_{1} \in \cP$ such that $P_{1} \leq P \leq P_{max}.$ Furthermore, there exist $d \in \NN$ and nonzero pairwise orthogonal
minimal $P_1,...,P_d \in \cP$ such that $P = \sum_{j=1}^{d} P_j.$
\end{cor}
\textbf{Proof.} The first assertion follows since otherwise there
exists an infinite sequence of projections that either strictly
increase or strictly decreases thus contradicting Lemma \ref{AN}.
The second assertion follows by recursively applying the first
assertion to the projection $P-P_1$ whenever $P-P_1 \neq 0.$
\hfill$\square$\ 
For any Hilbert space $H$ the set $\cB(H),$ consisting of bounded
operators on $H,$ forms a $C^*$--algebra under the adjoint and
operator norm. A $C^*$--subalgebra $\cA$ of $\cB(H)$ is called nondegenerate
if $\overline {\cA \, H} = H.$
\begin{lem}
\label{GNS}
Every $C^*$ algebra $\cA$ is isomorphic to a nondegenerate $C^*$--subalgebra
$\cA$ of $\cB(H)$ for some Hilbert space $H.$
\end{lem}
\textbf{Proof.} This isomorphism is given by the
Gelfand-Naimark-Segal construction \cite{murphy}.
\hfill$\square$\\

\medskip
Througout the remainder of this paper we identify $\cA$ with a
nondegenerate $C^*$--subalgebra of $\cB(H)$ for some Hilbert space $H$ as justified by Lemma \ref{GNS}. If $H_1$ and $H_2$ are closed subspaces of $H$ that are orthogonal we write $H_1 \perp H_2$ and if in addition their sum equals $H$ we write $H = H_1 \bigoplus H_2$ and $H_2 = H_{1}^{\perp}.$ We observe that for every $P, Q \in \cP,$ $PH$ is a closed subspace of $H,$ $P$ is the operator that orthogonally projects vectors onto $PH,$ $\hbox{ker}\, (P) = (PH)^{\perp},$ and $P \perp Q$ iff $PH \perp QH.$ We also observe that for every $A \in \cA_s,$
    $H = \hbox{ker}\, (A) \bigoplus \overline {AH}.$
\begin{cor}
\label{unital}
If $\hbox{spec}\, \, (A)$ is finite for every $A \in \cA_s$ then $\cA$ is unital.
\end{cor}
\textbf{Proof.} Lemma \ref{spec} and Corollary \ref{maxmin} ensure
that there exists a maximal projection $P \in \cP.$ It suffices to
show that $PH = H.$ Assume to the contrary that $PH \neq H.$ Since
$\cA$ is nondegenerate there exists $A \in \cA$ such that $AH$ is
not a subset of $PH$ and hence
$\hbox{ker} (P)$ is not a subset of $\hbox{ker} (A) = \hbox{ker} (A^*A).$
Therefore, since $Q \in \cA_s$ defined by
    $Q = P + A^*A$
satisfies
    $\hbox{ker} (Q) = \hbox{ker} (P) \bigcap \hbox{ker} (A),$
it follows that $\hbox{ker} (Q)$ is a proper subset of $\hbox{ker} (P)$
and hence
    $PH$
is a proper subspace of
    $\overline {QH}.$
Lemma \ref{spec} implies that $\overline {QH} = QH$ and that there
exists $T \in \cP$ such that $TH = QH.$ This implies that $P < T$
thus contradicting the assumption that $P$ is maximal.
\hfill$\square$\
\section{Derivation of Proposition \ref{pro2}}
\label{derivations2} \setcounter{equation}{0}
If $H_1$ and $H_2$ are closed subspaces of $H$ we let $Hom_{\cA}(H_1,H_2)$ denote the subspace consisting of $A \in \cA$ such that $AH_1 \subset H_2$ and such that the restriction of $A$ to $H_{1}^{\perp}$ equals $0.$ An element $U \in \cA$ is called a partial isometry if $U^*U \in \cP.$
\begin{lem}
\label{polar}
Every $A \in \cA$ can be factored uniquely as
$A = U|A|$
where
$H_1 = \left(\hbox{ker} (A)\right)^{\perp},$
$H_2 = \overline {AH},$
$|A| \in Hom_{\cA}(H_1,H_1)$ is the semipositive operator such that $|A|^2 = A^*A,$
and
$U \in Hom_{\cA}(H_1,H_2)$ is a surjective partial isometry such that $U^*U$ is orthogonal projection onto $H_1.$
\end{lem}
\textbf{Proof.} This is the polar decomposition theorem
\cite{murphy}.\hfill$\square$\\

\medskip Throughout this remainder of this section we will assume
that $\cA$ satisfies the property that $A \in \cA_s$ implies that
${spec}\, (A)$ is finite. Therefore Lemma \ref{spec} and Lemma
\ref{polar} ensures that $AV$ is closed whenever $V$ is a closed
subspace of $H$ and $A \in \cA.$
\begin{lem}
\label{dim1}
    If $P \in \cP$ is minimal and $A \in Hom_{\cA}(PH, PH)$ then there exists
    $c \in \CC$ such that $A = cP.$ Therefore, since $P$ is a nonzero element
    in $Hom_{\cA}(PH, PH),$ it follows that
$\dim \left( Hom_{\cA}(PH, PH) \right) = 1.$
\end{lem}
\textbf{Proof.}
    Since $A_r = \frac{1}{2}\left(A + A^*\right)$ and $A_i = \frac{1}{2i}\left(A - A^*\right)$ are in $\cA_s \cap Hom_{\cA}(PH, PH),$ Lemma \ref{spec} implies that $A_r$ and $A_i$ are either zero or linear combinations having nonzero real coefficients of projections. Each projection is in $Hom_{\cA}(PH, PH)$ and therefore is $\leq P.$ Since $P$ is minimal each projection must vanish or equal $P.$ Since $A = A_r + iA_i$ the proof is complete.
\hfill$\square$\
\begin{cor}
\label{dim2}
    If $P, Q \in \cP$ are minimal then
    $\dim \left( Hom_{\cA}(PH, QH) \right) \in \{0,1\}.$
\end{cor}
\textbf{Proof.}
    If nonzero elements $A, B \in Hom_{\cA}(PH,QH)$ have polar decompositions $A = U|A|$ and $B = O|B|$ then $O^*U \in Hom_{\cA}(PH,QH)$ and Lemma \ref{dim1} implies that $|A|,$ $|B|,$ and $O^*U$ are multiples of $P.$ This completes the proof.
\hfill$\square$\
\begin{pro}
\label{pro2}
$\cA$ is finite dimensional.
\end{pro}
\textbf{Proof.} Corollary \ref{unital} and Lemma \ref{maxmin}
imply that $I \in \cA$ and that there exists minimal projections
$P_1,...,P_d$ such that $I = \sum_{j=1}^{d} P_j.$ Therefore every
$A \in \cA$ can be expressed as the sum
$A = \sum_{i,j=1}^{d} A_{ij}$
where
$A_{ij} = P_i A P_j \in Hom_{\cA}(P_jH,P_iH).$
Therefore
$\cA = \bigoplus_{i,j = 1}^{d} Hom_{\cA}(P_jH,P_iH)$
and hence Corollary \ref{dim2} implies that $\dim \cA \leq d^2.$
\hfill$\square$\
\begin{rem}
\label{matrep}
In Corollary \ref{dim2} we can define an equivalence relation on the set $\{1,...,d\}$ by $i \sim j$ iff $\dim \left( Hom_{\cA}(P_jH,P_iH) \right) = 1.$ This partitions the set into $n \in \NN$ equivalence classes having cardinalities $d_1,...,d_n \in \NN$ with $d = \sum_{k = 1}^{n} d_k.$ Therefore $\cA$ is isomorphic to
$\bigoplus_{k = 1}^{n} \cB(\CC^d),$
which is a sum of matrix algebras, and
$\cA$ is abelian iff each $d_k = 1.$
\end{rem}
\section{Conclusions and Future Research}
\label{conclusions} \setcounter{equation}{0}
This paper derives simple criteria for a $C^*$--algebra to be
finite dimensional and describes the structure of the algebra
using properties of their projections that are related to fixed
point and spectral theory. A sharpened understanding of these
properties may be useful in quantum mechanics where
$C^*$--algebras provide the core mathematical foundation
\cite{strocchi}. Future research will examine the geometry of more complicated topological algebras that include the approximately finite dimensional $C^*$--algebras defined by Bratteli \cite{bratteli} and classified by Elliott \cite{elliott}, and the abelian Banach $*$--algebras whose fixed point properties were studied by
Dhompongsa and Fupinwong \cite{dhompongsa-fupinwong}.
\\ \\
{\it Acknowledgments} The authors thank Kazimierz Goebel for
helpful discussions.

\end{document}